\documentclass{amsart}
\usepackage{amssymb}
\usepackage{amsmath}
\usepackage{amsfonts}
\usepackage{pstricks}
\usepackage{pst-node}
\usepackage{graphicx}
\theoremstyle{plain}
\newtheorem{thm}{Theorem}[section]

\newtheorem{dfn}[thm]{Definition}

\newtheorem{qws}[thm]{Question}
\newtheorem{ass}[thm]{Assertion}
\newtheorem{prb}[thm]{Problem}
\newtheorem{cng}[thm]{Conjecture}
\theoremstyle{remark}

\def\pmc#1{\setbox0=\hbox{#1}
    \kern-.1em\copy0\kern-\wd0
    \kern.1em\copy0\kern-\wd0}

\begin{document}

\bigskip

\title[On generalized $3$-manifolds]{ On generalized $3$-manifolds which are not
homologically locally connected}

\bigskip

\author[U.~H.~Karimov]{Umed H. Karimov}
\address{Institute of Mathematics, 
Academy of Sciences of Tajikistan, 
Ul. Ainy $299^A$, Dushanbe 734063, 
Tajikistan}
\email{umedkarimov@gmail.com}

\author[D.~Repov\v{s}]{Du\v san Repov\v s}
\address{Faculty of Education,
and Faculty of  Mathematics and Physics,
University of Ljubljana, 
P.O.Box 2964, 
Ljubljana 1001, 
Slovenia}
\email{dusan.repovs@guest.arnes.si}


\begin{abstract}
We show that the classical example $X$ 
of a 
3-dimensional generalized manifold
constructed by van Kampen is another example of not
homologically locally connected (i.e. not HLC) space.
This space $X$ is
not locally homeomorphic to any of
the compact metrizable 3-dimensional
manifolds constructed in our earlier paper  
which are not HLC
spaces either.
\end{abstract}

\date{\today}

\maketitle

\section{Introduction}

In our earlier paper \cite{KR}  we constructed for every natural number $n > 2$, examples
of $n$-dimensional compact metrizable cohomology $n$-manifolds
which are not homologically locally connected with respect to
the
singular homology  (i.e. they are not  {\it HLC} spaces). In the present paper we
shall call them {\it singular quotient $n$-manifolds.}

Subsequently, we have discovered that van Kampen constructed a
compact metrizable generalized $3$-manifold which "is not locally
connected in dimension 1 in the {\sl homotopy sense}" \cite[p.
573]{Bg}. The description of van Kampen's construction can be
found in \cite[p. 573]{Bg}
(see also \cite[p. 245]{W}).

The obvious modification of van Kampen's construction gives an
infinite class of examples -- we shall call them {\it van Kampen
generalized $3$-manifolds}. The main purpose of the present paper
is to prove the following theorem:

\begin{thm}\label{Theorem 1}
No van Kampen generalized $3$-manifold
is homologically locally
connected in dimension 1 with respect to
the singular homology. They
are neither locally homeomorphic to any singular quotient
$3$-manifold. Furthermore, no singular quotient $3$-manifold is
locally homeomorphic to any van Kampen generalized $3$-manifold.
\end{thm}

\section{Preliminaries}

We shall denote the singular homology
(resp.  \v{C}ech cohomology)
groups with integer coefficients
by $H_{\ast}$ (resp.  $\check{H}^{\ast}$). 
All spaces considered in this paper will be assumed
to be metrizable, locally compact and finite-dimensional.   Under
these circumstances the classes of {\it generalized manifolds},
{\it homology manifolds}, and {\it cohomology manifolds} in the classical
sense
\cite{W} coincide (cf. \cite{Br, Hr, Mit}). For the
history and importance of generalized manifolds see \cite{R2}.

\begin{dfn} {\rm(cf. \cite[p. 377, Corollary 16.9]{Br}).}\label{Generalized manifold}
A locally compact, cohomologically locally connected with respect
to \v Cech cohomology $(clc)$,and
cohomologically finite-dimensional
space $X$ is called a {\sl generalized $n$-manifold},
$n \in \mathbb N$,
if
$$\check{H}^p(X, X \setminus \{x\}) \cong
\check{H}^p(\mathbb{R}^n, \mathbb{R}^n \setminus \{0\})$$ for all
$x \in X$ and all $p \in \mathbb{Z}_{+}$.
\end{dfn}

\begin{dfn}\label{Homology 3-sphere}
A closed 3-manifold is called a {\sl homology $3$-sphere} if its homology
groups are isomorphic to the homology groups of the standard
$3$-sphere {\rm (}cf. e.g. \cite{G}{\rm )}. The complement of the
interior of any 3-simplex in any triangulation of any homology
$3$-sphere is called a {\sl homology $3$-ball}.
\end{dfn}

In 1904 Poincar\'e constructed 
the first example of a homology $3$-sphere
with a nontrivial fundamental group  \cite{P, ST}. The
following theorem  is due to Brown \cite{Bn}:

\begin{thm}\textrm{\bf{(Generalized Sch\" onflies theorem)}}. \label{Schoenflies}
Let $h:S^{n-1}\times [-1,1]\to S^n$ be an embedding. Then the
closures of both complementary domains of $h(S^{n-1}\times 0)$ in
$S^n$ are topological $n$-cells.
\end{thm}

\begin{dfn}\label{Local homeomorphism}
A space $X$ is said to be {\sl locally homeomorhic} to the space
$Y$ if for every point $x\in X$ there exists an open neighborhood
$U_x \subset X$ which is homeomorphic to some open subspace of \
$Y. $
\end{dfn}

\begin{dfn}\label{Commutator length}{\rm(cf. e.g.
\cite{Ba, EKR, KR})}. Let $G$ be any group and  $g\in G$ any element.
The {\sl commutator length} \ $cl(g)$ of the element
$g$ is defined as the minimal number of the commutators of the
group $G$ the product of which is $g.$ If such a number does not
exist then we set $cl(g) = \infty.$ Moreover, $cl(g) = 0$ if and
only if $g = e$, {where $e$ denotes the neutral element}.
\end{dfn}

\begin{thm}\label{Griffiths}{\rm(cf. \cite{GT, Gr}).}
Let $G$ be a free product of groups $\{G_i\}_{i=1}^k,$ \ $G = G_1*
G_2 *\cdots *G_k.$ Let $g\in G,$ $g_i\in G_i$ for $i =
\overline{1,k}$ and let $g = g_1g_2\cdots g_k$. Then $cl(g) =
\sum_{i=1}^k cl(g_i).$
\end{thm}

\begin{ass}\label{Property 2 of cl}
The commutator length function has following properties:

$($1$)$ If $\varphi: G \to H$ is a homomorphism of groups and $g
\in G$ then $cl(\varphi(g)) \leq cl(g);$ and

$($2$)$ For every element $g\in G, \ cl(g) = \infty$
$\Leftrightarrow g \notin [G, G].$
\end{ass}

All undefined terms can be found in \cite{Br}, \cite{ST} or
\cite{S}.

\bigskip

\section{The construction of the
van Kampen example and  the proof of the Main Theorem}

To prove the main theorem we shall need a more detailed
description of the original van Kampen construction \cite[p.
573]{Bg}. For nonnegative numbers $i,p$ and $q$ let the {\it ball
layer} $A_{i,p,q}$ be the following subspace of $\mathbb R^3$:
$$A_{i,p,q} = \{\overline{a}\in \mathbb{R}^3 \ | \ p \leq
|\overline{a} - (0, i, 0)| \leq p + q, \ (0,i,0) \in \mathbb{R}^3
\}.$$ In particular, $A_{i,0,q}$ is a 3-ball with the center at
the point $(0,i,0)\in \mathbb{R}^3$ and of radius $q.$ For $i\in
\mathbb{N}$, let $\{B_i\}_{i \in \mathbb N}$ be a sequence of a homology $3$-balls
with a nontrivial fundamental group. All spaces $B_i$ are compact
$3$-manifolds the boundaries of which are homeomorphic to $S^2.$

Let for $n,  k \in \mathbb{N}$
$$X_{n, k} =
\overline{A_{0,(n - \frac{1}{2}), k}\setminus \bigcup
_{j=n}^{n+k-1} A_{j, 0, \frac{1}{4} }}\ \bigcup \
\bigcup_{i=n}^{n+k-1}B_i$$ with the topology of identification of
the boundaries of $A_{i, 0, \frac{1}{4} }$ and $B_i$ \ (the 
bar over the space denotes the closure of
that space). Let
$$X_k = X_{1, k}\cup A_{0,0,\frac{1}{2}}.$$
Obviously, $X_{n, k}$ and $X_k$ are compact $3$-manifolds with
boundaries. Since the homology groups of $B_i$ are the same as the
homology groups of the 3-ball $A_{i, 0, \frac{1}{4} }$ it follows
from the exactness of the Mayer-Vietoris homology sequence that
the homology groups of $X_k$ are the same as the homology groups
of $A_{0,0, (k + \frac{1}{2})},$ i.e. they are trivial and $X_k$
is a
homology $3$-ball.


Consider $X = \underrightarrow{\textrm{lim}} (X_k \subset
X_{k+1})$ with the topology of the direct limit. Then
\begin{equation}\label{Acyclicity}
H_*(X) \cong H_*(*) \ \
\hbox{\rm and} \ \
\check{H}^*(X) \cong \check{H}^*(*), \
\text{i.e.} \ X \ \text{is an acyclic $3$-manifold}.
\end{equation}

The one-point compactification of this space, by the point $\ast$,
is a
3-dimensional compact metrizable  space $X^{\ast}$.

We have
\begin{equation}\label{X^{ast}}
X^{\ast} \cong \underleftarrow{\textrm{lim}}(X_k/\partial A_{0,0,
(k + \frac{1}{2})}),
\end{equation}
where $X_k/\partial A_{0,0, (k + \frac{1}{2})}$ is the quotient of
the manifold $X_k$ via the identification of all of its boundary
points to one point and where the projections $$X_k/\partial
A_{0,0, (k + \frac{1}{2})}\leftarrow X_{k+1}/\partial A_{0,0, (k +
1 + \frac{1}{2})}$$ map the subspaces $\overline{(X_{k+1}/\partial
A_{0,0, (k + 1 + \frac{1}{2})})\setminus X_k}$ of the spaces
$X_{k+1}/\partial A_{0,0, (k + 1 + \frac{1}{2})}$ to the point of
identification of $X_k/\partial A_{0,0, (k + \frac{1}{2})}$.

Obviously, the spaces $X_k/\partial A_{0,0, (k+ \frac{1}{2})}$ are
homology $3$-spheres and all projections in the inverse sequence
generate isomorphisms of homology and cohomology groups. Therefore
\begin{equation}\label{Homol sphere}
\check{H}^*(X^{\ast}) \cong \check{H}^*(S^3).
\end{equation}

It follows from the exact sequence
$$\cdots\rightarrow \check{H}^{p-1}(X^{\ast}\setminus \ast)
\rightarrow \check{H}^{p}(X^{\ast}, X^{\ast}\setminus
\ast)\rightarrow \check{H}^p(X^{\ast})\rightarrow \check{H}^
p(X^{\ast}\setminus \ast)\rightarrow \cdots,$$ \noindent from the
facts that $X^{\ast}\setminus \ast = X, \ \check{H}^p(X^{\ast}) =
0 \ \text{for} \ p \neq 3, $ and from (\ref{Acyclicity}) and
(\ref{Homol sphere}), that
\begin{equation}\label{Local groups}
 \check{H}^{p}(X^{\ast}, X^{\ast}\setminus \ast) \cong
H^p(\mathbb{R}^3, \mathbb{R}^3 \setminus \{0\}).
\end{equation}

Let us show that $X^{\ast}$ is a $clc$  space. Since $X$ is a
$3$-manifold, the space $X^{\ast}$ is a $clc$ space at every point
of $X^{\ast}\setminus \ast.$ Since $X =
\underrightarrow{\textrm{lim}} (X_k \subset X_{k+1})$ and $X_k$
are compact spaces, the system of the sets $\{X^{\ast}\setminus
X_k\}$ is a basis
of neighborhoods of the point $\ast.$ According
to (\ref{X^{ast}}), we have
$$\overline{X^{\ast}\setminus
X_k} = \underleftarrow{\textrm{lim}}\ (\overline{X_n/\partial
A_{0,0, (n + \frac{1}{2})})\setminus X_k}\ , \ \textrm{for} \ \ n
> k. $$

\noindent The spaces $\overline{(X_n/\partial A_{0,0, (n +
\frac{1}{2})})\setminus X_k}$ are homology $3$-balls. It follows
that $\overline{X^{\ast}\setminus X_k}$ is an acyclic space with
respect to the \v{C}ech cohomology. Therefore $X^{\ast}$ is a
$clc$ space.
Hence it follows from (\ref{Local groups}) that $X^{\ast}$ is
a generalized $3$-manifold (cf. Definition \ref{Generalized
manifold}). 

The topological type of the space $X^{\ast}$ depends
on the choice of the sequence of homology $3$-balls $B_i$. In the
case when all $B_i$ are homeomorphic to the Poincar\'e homology
$3$-sphere, the space $X^*$ is the example of van Kampen. In
general, since there exist infinitely many distinct homology
$3$-balls (cf. \cite{G}),
there exists an infinite class of van Kampen
generalized $3$-manifolds.

Let $X^{\ast}$ be any van Kampen generalized $3$-manifold. We shall
show that $X^{\ast}$ is not $HLC$ in dimension $1$. It suffices to
prove that for every $k$ the embedding $X^{\ast}\setminus X_k
\subset X^{\ast}\setminus X_1$ is not homologically trivial with
respect to singular homology.

The space $\overline{X^{\ast}\setminus X_k }$ is a retract of
$\overline{X^{\ast}\setminus X_1}.$ Indeed, the $3$-balls $A_{i,0,
\frac{1}{4}}$ are $AR$ spaces so there exists for every $i$, a
mapping of $B_i$ on $A_{i,0, \frac{1}{4}}$ which is the identity
on the boundary of $B_i$. Therefore we have a mapping of $X_{1,k}$
onto the ball layer $A_{0, \frac{1}{2}, k}$.

However, the sphere $A_{0, k+\frac{1}{2}, 0 }$ is a retract of
this ball layer and $A_{0, k+\frac{1}{2}, 0 }\subset
\overline{X^{\ast}\setminus X_k} $. Therefore
$\overline{X^{\ast}\setminus X_k} $ is a retract of
$X^{\ast}\setminus X_1$. So in order to prove that
$\overline{X^{\ast}\setminus X_k }\subset X^{\ast}\setminus X_1$
is not homologically trivial it suffices to prove that
$H_1(\overline{X^{\ast}\setminus X_k}) \neq 0.$

Since for any path-connected space the $1$-dimensional singular
homology group is isomorphic to the abelianization of the
fundamental group of this space, it
suffices
to show that
$\pi_1(\overline{X^{\ast}\setminus X_k})$ is a nonperfect group.

Consider the union of spheres $A_{j, \frac{1}{4}, 0},$ for $j > k,$
with the compactification point as a subspace of the space
$X^{\ast}\setminus X_k.$ Join them by the segments
$$\rm{I}_j =
\{(0, y, 0)|\ y\in [j + \frac{1}{4}, j + 1 - \frac{1}{4}]\},  
\ \ 
j
= k+1, \cdots \infty .$$
We get the compactum $$A = \bigcup_ {j = k+1}^{\infty}A_{j,
\frac{1}{4}, 0}\cup  \rm{I}_j\cup\{\ast\}.$$

The space $A$ also lies in a $3$-dimensional cube with a countable
number of open balls $\{\textrm{int}A_{j, 0, \frac{1}{4}}\}_{j =
k+1}^{\infty}$ removed and it is obviously its retract. Therefore
the space
$$B = A\bigcup\bigcup_{i=k + 1}^{\infty}B_i $$
\noindent with the natural topology of identification is a retract
of $X^{\ast}\setminus X_k.$

Consider nontrivial loops $\alpha_i \in B_i,$ with their homotopy
classes $[\alpha_i].$  Then $0 < l([\alpha_i]).$ Since $H_1(B_i) =
0$ it follows that $l([\alpha_i]) < \infty.$ Let $\alpha$ be any
loop in $B$ such that for the canonical projection of this space
to $B_i$, the image of $\alpha$ generates the loop $\alpha_i$
(obviously such a
loop exists). We have a homomorphism 
$$\pi_1(B) \to
\pi_1(B_1)*\pi_1(B_2)*\cdots$$
such that the image of $[\alpha]$
has projections $[\alpha_i]$ in $\pi_1(B_i).$
Hence
according
to 
\ref{Griffiths}
and 
\ref{Property 2 of cl}, 
the commutator
length of $[\alpha]\in \pi_1(B)$ is infinite, $cl([\alpha]) =
\infty$, and the group $\pi_1(B)$ is not perfect (\ref{Property 2
of cl}). Therefore and since $B$ is a retract of
$X^{\ast}\setminus X_k$ it follows that
$\pi_1(\overline{X^{\ast}\setminus X_k})$ is nonperfect and hence
the van Kampen generalized $3$-manifold is not an
$HLC$-space.

Let us prove that no van Kampen generalized $3$-manifold is
locally homeomorphic to any singular quotient $3$-manifold.
Consider any open neighborhood $U$ of the singular point $\ast$ of
a van Kampen generalized $3$-manifold. Suppose that it were
homeomorphic to some open subset of some almost $3$-manifold. Now,
every singular quotient $3$-manifold is a quotient space of a
topological $3$-manifold by some continuum, which generates a
singular point \cite{KR}. Therefore there should exists an index
$i$ such that some neighborhood of $B_i $ in $U$ embeds into the
Euclidean $3$-space with some ball layer $A_{i,\frac{1}{4},
\varepsilon}$.

By Theorem
\ref{Schoenflies}, the bounded component of $\mathbb{R}^3\setminus
A_{i, \frac{1}{4}+\frac{\varepsilon}{2}, 0}$ is a $3$-cell. So the
space $B_i\cup A_{i,\frac{1}{4}+\frac{\varepsilon}{2}, 0}$ embeds
in a $3$-cell. This embedding cannot be an onto mapping to the
$3$-cell because $B_i$ is not simply connected. Therefore the
sphere $A_{i, \frac{1}{4}+\frac{\varepsilon}{2}, 0}$ must be a
retract of $B_i$. However, this is impossible since the space
$B_i$ is acyclic whereas the sphere $A_{i, \frac{1}{4}, 0}$ is not
acyclic. So no van Kampen generalized $3$-manifold
is locally
homeomorphic to any almost $3$-manifold.

By the argument above and since every local homeomorphism must map
singular points to singular points it follows that no singular
quotient $3$-manifold is locally homeomorphic to any van Kampen
generalized $3$-manifold. \qed

\section{Epilogue}

Cohomological local connectedness does not imply local
contractibility even
in the category of compact metrizable
generalized $3$-manifolds (this follows e.g. from
Theorem~\ref{Theorem 1}). However, the following problem remains
open:
\begin{prb}
Does there exist a finite-dimensional compact metrizable
generalized manifold (or merely finite-dimensional compactum)
which is homologically locally connected but not locally
contractible?
\end{prb}

Marde\v{s}i\'{c} formulated the following interesting problem:

\begin{qws}
Is it true that every compact generalized $n$-manifold which is
an
ANR
can be represented as an inverse limit of compact $n$-manifolds?
\end{qws}

As it was mentioned above (Chapter 3 of reference (\ref{X^{ast}})),
every van Kampen generalized $3$-manifold is an inverse limit of
compact $3$-manifolds. In contrast with this we set forth
the following 
conjecture:
\begin{cng}
No singular quotient $3$-manifold can be represented as an inverse
limit of closed $3$-manifolds.
\end{cng}

\section{Acknowledgements}
This research was supported by the Slovenian Research Agency
grants P1-0292-0101, J1-2057-0101 and J1-4144-0101.

\end{document}